\documentclass[11pt]{article}
\catcode`\@=11 \catcode`\@=12

\usepackage[utf8]{inputenc}
\usepackage{oldgerm}
\usepackage{amssymb,amsmath,amsthm,stmaryrd}
\usepackage{mathrsfs}
\usepackage{graphicx,psfrag,epsfig}
\usepackage[pdftex,colorlinks,urlcolor=blue,pdfstartview=FitH]{hyperref}
\usepackage[all]{xy}
\usepackage{pst-all}
\usepackage{multicol}
\usepackage{float}
\usepackage{enumerate}
\newcommand{\Rm}{\mathbb{R}}

\newcommand{\Qm}{\ensuremath{\mathbb{Q}}}
\newcommand{\Nm}{\ensuremath{\mathbb{N}}}
\newcommand{\Zm}{\ensuremath{\mathbb{Z}}}

\newcommand{\mM}{\ensuremath{\mathcal{M}}}

\newcommand{\mA}{\ensuremath{\mathcal{A}}}

\newcommand{\Tm}{\ensuremath{\mathbb{T}}}
\newcommand{\vs}{\vspace{.2cm}}
\def\proof {\noindent{\sc{Proof. }}}
\def\qed {\mbox{}\hfill {\small \fbox{}} \\}
\def\lto{\longrightarrow}

\def\leq{\leqslant}
\def\geq{\geqslant}

\newcommand{\T}{\text{\textbf{T}}}

\newtheorem{thm}{Theorem}

\newtheorem{prop}{Proposition}[section]
\newtheorem{lem}{Lemma}[section]

\newtheorem{cor}{Corollary}[section]
\newtheorem{pte}{Property}[section]

\theoremstyle{definition}
\newtheorem{Def}{Definition}[section]

\theoremstyle{remark}
\newtheorem{rem}{Remark}[section]
\newtheoremstyle{restate}
{\topsep}
{\topsep}
{\itshape}
{}
{\bfseries}
{.}
{ }
{\thmname{#1}\thmnote{ #3}}
\theoremstyle{restate}

\newtheorem{thm*}{Theorem}
\newtheorem{prop*}{Proposition}
\newtheorem{lem*}{Lemma}
\newtheorem{add*}{Addendum}
\newtheorem{cor*}{Corollary}
\newtheorem{pte*}{Property}
\theoremstyle{definition}
\newtheorem{Def*}{Definition}
\newtheorem{exm*}{Example}
\theoremstyle{remark}
\newtheorem{rem*}{Remark}

\theoremstyle{definition}

\theoremstyle{remark}

\newcommand{\N}{\mathbb N}

\newcommand{\R}{\mathbb R}

\newcommand{\setm}{\setminus}

\newcommand{\eps}{\varepsilon}

\newcommand{\vp}{\varphi}

\newcommand{\Sig}{\Sigma}

\newcommand{\Tg}{\text{\textbf{T}}}

\newcommand{\rit}{\rightarrow}
\newcommand{\ma}{\mapsto}

\newcommand{\lmto}{\longmapsto}

\newcommand{\hto}{{\rm{h_{top}}}}
\newcommand{\hp}{{\rm{h_{pol}}}}

\newcounter{paraga}[subsubsection]

\author{Patrick Bernard,   Clémence Labrousse}

\title{An entropic characterization of the flat metrics on the two torus}

\begin{document}

\maketitle

\begin{center}
-----
\end{center}
\begin{small}
\begin{multicols}{2}

\noindent
Patrick Bernard
\footnote{Universit\'e  Paris-Dauphine},\\
\'Ecole Normale sup\'erieure,\\
DMA (UMR CNRS 8553)\\
45, rue d'Ulm\\
75230 Paris Cedex 05,
France\\
\texttt{patrick.bernard@ens.fr}\\

\noindent
Clémence Labrousse
\footnotemark[1],\\
\'Ecole Normale sup\'erieure,\\
DMA (UMR CNRS 8553)\\
45, rue d'Ulm\\
75230 Paris Cedex 05,
France\\
\texttt{clemence.labrousse@ens.fr}\\

\end{multicols}
\vs
\thispagestyle{empty}
\begin{center}
-----
\end{center}

\textbf{Abstract. } 
The geodesic flow of the flat metric on a torus is minimizing the polynomial entropy among 
all geodesic flows on this torus. We prove here that this properties characterises the flat metric
on the two torus.

\begin{center}
-----
\end{center}

\textbf{R\'esum\'e. }
Le flot géodésique des métriques plates sur un tore minimise l'entropie polynomiale parmi tous les flots
géodésiques sur ce tore.
On montre ici que cette propriété caractérise les métriques plates en dimension deux.

\begin{center}
-----
\end{center}

\tableofcontents

\vfill
\hrule
\vs
{The research leading to these results has received funding from the European Research Council
under the European Union's Seventh Framework Programme (FP/2007-2013) / ERC Grant
Agreement  307062. }

\end{small}

\newpage

\section{Introduction}

There are several classes of hyperbolic Manifolds on which 
the metrics with constant curvature are characterized 
by the fact that their geodesic flow is minimizing the 
topological entropy, see \cite{Kat} and \cite{BCG-96} for example.
The situation is different on tori.
Flat metrics have zero entropy, but other metrics also have zero 
entropy, such as the tori of revolution.
In order to characterize the flat metrics, it is therefore
useful to consider a finer dynamical invariant of the geodesic
flow, such as the polynomial entropy, introduced in \cite{Mar13}.

Using the techniques of  \cite{Mar13}, it was proved in \cite{L12a}
that the polynomial entropy of a flat 
torus of dimension $d$ (in restriction to the sphere bundle)
is equal to $d-1$, which is a lower bound for the polynomial 
entropy of all metrics on $\Tm^d$.
It was also proved in \cite{LM} that the polynomial entropy 
of the revolution two torus is two, which is higher than the one
of the flat two tori.
This gives an indication that the polynomial entropy might be a sufficiently fine
invariant to characterize the flat metric.
Our main result in the present paper is that this is indeed the case in dimension two.
A partial result in that direction has been obtained in \cite{L12b}.

\begin{thm}\label{maingeo}
If the polynomial entropy of a $C^2$ metric $g$ on $\Tm^2$
(in restriction to the sphere bundle) is smaller than two,
then this entropy is equal to one and the torus $(\Tm^2, g)$ is isometric to a flat torus.
\end{thm}

We will prove this result using Mather-Fathi theory.
The useful facts from this theory are recalled in Section \ref{sectonelli},
where a more general estimate on the polynomial entropy of Tonelli Hamiltonians 
is given, see Theorem \ref{maintonelli}.
Theorem \ref{maingeo} is deduced from  Theorem \ref{maintonelli}
using the theorem of Hopf and its variants, see  \cite{Innami}.
The definition of the polynomial entropy is recalled in 
Section \ref{secentro}, and the entropy estimates leading 
to the proof of  Theorem \ref{maintonelli} are detailed in Section \ref{secproof}.
Once the dynamics has been well understood with the help of Mather-Fathy theory,
these estimates are similar to those appearing in \cite{Mar13}, \cite{LM}, and \cite{Mar14}.

\section{The polynomial entropy}\label{secentro}

Consider a continuous  map $f:X\rit X$, where $(X,d)$ is a compact metric space.
We construct new metrics $d_n^f$ on $X$
 by setting 
\[
d_n^f(x,y)=\max_{0\leq k\leq n-1}d(f^k(x),f^k(y)).
\]
These metrics are the \textit{dynamical metrics} associated with $f$. 
Obviously, if $f$ is an isometry or is contracting, $d^f_n$ coincides with $d$ and 
in general $d_n^f$ is topologically equivalent to $d$.
We denote by $G_n^f(\eps)$ the minimal number of balls of radius $\eps$ for the metric $d_n^f$ in a finite covering of $X$. 
The topological entropy of the map $f$, defined as,
$$
\hto(f)=\lim_{\eps\rit 0}\limsup_{n\rit\infty}\frac{\log G_n^f(\eps)}{ n}
$$
measures the exponential growth rate of $G_n^f$.   In the 
present paper we will rather consider a polynomial measure
of the growth rate introduced in \cite{Mar13}:

\begin{Def}
The \emph{polynomial entropy} $\hp(f)$ of $f$  is defined by
\begin{align*}
\hp(f) =\lim_{\eps\rit 0}\limsup_{n\rit\infty}\frac{\log G_n^f(\eps)}{\log n}.
\end{align*}
\end{Def}

We also consider sets that are $\eps$\emph{-separated}  for the metrics $d_n^f$ (we will write $(n,\eps)$-separated). Recall that a set $E$ is said to be $\eps$-separated for a metric $d$ if  for all $(x,y)$ in $E^2$, $d(x,y)\ge \eps$. 
Denote by  $S_n^f(\eps)$ the maximal cardinal of a $(n,\eps)$-separeted set contained  in $X$. 
Observing that 
$
S_n^f(2\eps)\leq G_n^f(\eps)\leq S_n^f(\eps),
$
we obtain
\[
\hp(f) =\lim_{\eps\rit 0}\limsup_{n\rit\infty}\frac{\log S_n^f(\eps)}{\log n}.
\]

\begin{rem}
If $\phi:=(\phi^t)_{t\in \R}$ is a continuous flow on $X$, for $t>0$ and $\eps>0$,  one can define in the same way the numbers $G_t^\phi(\eps)$ and $S_t^\phi(\eps)$.
The polynomial entropy  $\hp(\phi)$ of $\phi$  is defined as 
\[
\hp(\phi) =\lim_{\eps\rit 0}\limsup_{t\rit\infty}\frac{\log G_t^\phi(\eps)}{\log t}=\lim_{\eps\rit 0}\limsup_{t\rit\infty}\frac{\log S_t^\phi(\eps)}{\log t}.
\] 
One easily checks that if $\phi^1$ is the time-one map of $\phi$, $\hp(\phi)=\hp(\phi^1)$.
\end{rem}

The following properties of the polynomial entropy are proved in  \cite{Mar13}. 
\begin{pte}\label{ptehphw}

\begin{enumerate} 
\item $\hp$ is a $C^0$ conjugacy invariant, and does not depend on the choice of topologically equivalent
metrics on $X$.
\item If $A$ is a subset of $X$ $f$--invariant, $\hp(f_{\vert_A})\leq \hp(f)$.
\item For $m\in\N^*$, $\hp(f^m) =\hp(f)$ and if $f$ is invertible, $\hp(f^{-m})=\hp(f)$.
\end{enumerate}
\end{pte}

We conclude this section with the following 
useful result which relates the polynomial entropy of a flow with that of Poincar\'e map.

\begin{prop}\label{Poincare}
Let $M$ be a smooth manifold, $d$ a distance on $M$ associated with a Riemannian metric, and $X$ a 
$C^1$ complete vector field on $M$ with flow $\phi=(\phi_t)_{t\in \R}$. Let $A$ a be compact $\phi$--invariant subset of $M$ and let $\Sig$ be a $C^1$ codimension $1$ submanifold of $X$ such that: 
\begin{itemize}
\item for any $a\in A$, there exists $t>0$ such that $\phi^{t}(a)\in \Sig$. 
\item for any $a\in A\cap \Sig$, $X(a)$ is transverse to $\Sig$.
\end{itemize}
Then the Poincar\'e return map $\vp: A\cap\Sig \rit A\cap \Sig$ is well defined, continuous and satisfies
\[
\hp(\vp)\leq \hp(\phi,A).
\]
\end{prop}

\proof
 Let $\tau :A\cap \Sig \rit \R^*_+ :a \ma \tau_a$ be the first return time map of $\vp$.
  
Since the function $\tau$ is continuous 
on the compact set $A\cap \Sig$, we have 
$T:=\max\{\tau_a\,|\, a\in A\cap \Sig\}<\infty$. Let  $d_\Sig$ be the distance induced by $d$ on $\Sig$.

There exists $\tau^*>0$ and a neighborhood $V$ of $A\cap \Sigma$ in $\Sigma$
such that 
 the map $\Phi : ]-4\tau^*,4\tau^*[\times V\cap\Sig \rit M : (t,a)\ma \phi^t(a)$ is a 
 $C^1$--diffeomorphism onto its image.
 Its inverse is thus locally Lipschitz, hence its restriction to the compact set
  $K:=\Phi([-\tau^*,\tau^*]\times (A\cap\Sig))$ is Lipschitz.
 As a consequence, there exists $\delta>0$ such that 
  \begin{equation}\label{Lip}
d(x,x')\geq \delta \max (|t-t'|,d_\Sig(a,a')) \quad {\rm{if}}\quad x=\phi^t(a) \quad {\rm{and}}\quad x'= \phi^{t'}(a').
\end{equation}
Note that  $\tau^*<\frac{1}{4}\min\{\tau_a\,|\, a\in A \cap \Sig\}$.
Since the compact sets $A\cap \Sigma$ and 
$A\setminus\Phi(]-\tau^*,\tau^*[\times V)$ are disjoint, the constant $\delta$ can be chosen
such that
\begin{equation}\label{>eps}
d(a,x)\geq \delta \tau^* \quad \text{for each}\quad 
a\in \Sigma \cap A \text{ and }x\in A\setminus K.
\end{equation}
Let $\tau_x^k$ be the successive return times of the point $x$, so that $\varphi^k(x)=\phi^{\tau^k_x}(x)$.
Note that $\tau^1_x=\tau_x$, and $\tau^{k+1}_x=\tau^{k}_x+\tau_{\varphi^{k}(x)}$,
hence $\tau^k_x\leq kT$ for all $x\in A\cap \Sigma$.

We will now prove that two points $x$ and $y$ of $A\cap \Sig$ which are $(n,\eps)$--separated
by $\varphi$ are $(nT,\delta\epsilon)$ separated by $\phi$ for $\epsilon<\tau^*$.
There exists $m\in \{0,\dots,n\}$ such that $d_\Sig(\vp^m(x),\vp^m(y))\geq \eps$.
Let us assume for definiteness that $\tau^m_x\leq \tau^m_y$.

If $\phi^{\tau_x^m}(y)\in A\setm K$, 
then $d(\phi^{\tau^m_x}(x),\phi^{\tau^m_x}(y))\geq \delta\epsilon$
 by (\ref{>eps}), hence $x$ and $y$ are $(\tau^m_x,\delta \epsilon)$-separated
 by $\phi$.

If $\phi^{\tau_x^m}(y)\in  K$, then there exists $m'\leq m$ 
 and $s\in [-\tau^*,\tau^*]$ such that $\phi^{\tau^m_x}(y)=\phi^s(\vp^{m'}(y))$.

If $m'=m$, then $d(\phi^{\tau^m_x}(x),\phi^{\tau^m_x}(y))=d(\vp^m(x),\phi^{s}(\vp^m(y))\geq \delta \max(s,\eps)$, hence $x$ and $y$ are $(\tau^m_x,\delta \epsilon)$-separated
 by $\phi$.

If $m'<m$, then there exists  $k\in \{1,\dots,m\}$ such that $\phi^{\tau^k_x}(y)\notin K$,
which implies that $d(\phi^{\tau^k_x}(x), \phi^{\tau^k_x}(y))\geq \delta \eps$ by (\ref{>eps}).
As a consequence, the point $x$ and $y$ are $(\tau^k_x,\delta\eps)$-separated by $\phi$.

We have proved that 
$$
S^{\varphi}_n(\epsilon)\leq S^{\phi}_{nT}(\delta\epsilon)
$$
provided $\epsilon<\tau^*$, which implies the inequality on $\hp$.
\qed

\section{Tonelli Hamiltonians}\label{sectonelli}

\subsection{Some definitions from weak KAM theory}
We work on the $d$-dimensional  torus $\T:=\Rm^d/\Zm^d$, and will mostly consider the case $d=2$.
A Tonelli Hamiltonian on $\T$ is a $C^2$ Hamiltonian function 
$H(q,p): \T\times \Rm^d (=T^*\T)\lto \Rm$ such that, for each $q\in \T$, the function
$p\lmto H(q,p)$ is convex with positive definite Hessian and superlinear.
The Hamiltonian vectorfield on $T^*\T$ is given by
$$
X_H(q,p)=\big( -\partial_qH(q,p),\partial_pH(q,p)
\big).
$$
It generates a complete flow $\varphi_H^t$ which preserves the function $H$.
To each Tonelli Hamiltonian is associated 
the Lagrangian function $L$ on $T\T=\T\times \Rm^d$ given by
$$
L(q,v)=\sup _{p\in \Rm^d} \big( p\cdot v-H(q,p)\big)
$$
and the \textit{Legendre diffeomorphisms}
\begin{align*}
\T\times \Rm^d=T^*\T \ni (q,p)&\lmto (q,v)=(q,\partial_pH(q,p))\in T\T=\T\times \Rm^d,\\
\T\times \Rm^d=T\T \ni (q,v)&\lmto (q,p)=(q,\partial_vL(q,v))\in T^*\T=\T\times \Rm^d,
\end{align*}
which are inverse of each other.
To a Riemaniann metric
$g_x(v,v')=\langle G(x) v,v'\rangle$
where $G(x)$ is a $C^2$ field of positive definite symmetric matrices,
we associate the pair
$$
L(x,v)=\frac{1}{2}\langle G(x) v,v\rangle \quad, \quad
H(x,p)=\frac{1}{2}\langle G(x)^{-1} p,p\rangle.
$$
It is well-known that the Hamiltonian flow of $H$ is conjugated to the geodesic
flow by the Legendre diffeomorphism $(x,v)\lmto (x,G(x)v)$. 

Returning to the general case of a Tonelli Hamiltonian, 
the $\alpha$ function of Mather is defined on $H^1(\T,\Rm)$ by
$$
\alpha(c):= \inf_{u\in C^{\infty}} \sup _q H(q,c+du(q))
=\min_{u\in C^{1,1}}\sup _q H(q,c+du(q)),
$$
where the infimum and the minimum are taken respectively 
on the set of smooth functions on $\T$
and on the set of $C^1$ functions with Lipschitz differential.
It was proved in \cite{ENS} that the minimum exists 
on the set of $C^{1,1}$ functions, see also \cite{FS:04}.
A $C^{1,1}$ function  satisfying the inequality
$$
H(q,c+du(q))\leq \alpha(c)
$$
at each point is called a $c$-critical subsolution (as we just recalled, such functions
exist).
There may exist several $c$-critical subsolutions.
At least one of them, $w$, has the property that 
$$
H(q,c+dw(q))=\alpha(c) \Rightarrow H(q,c+du(q))=\alpha(c)
$$
for all  critical  subsolutions $u$. We define 
$$
\mA(c):=\{q\in \T, H(q,c+dw(q))=\alpha(c)\}
=\cap_u \{q\in \T, H(q,c+du(q))=\alpha(c)\},
$$
where the intersection is taken on all $c$-critical subsolutions $u$.
This is a non-empty compact set, called the projected Aubry set.
In view of the strict convexity of $H$ in $p$, the differential
$du(q)$ of a $c$-subsolution $u$ at a point $q\in \mA(c)$ does not 
depend on  the $c$-critical subsolution $u$.
We define
$$
\mA^*(c):=\{ (q,c+dw(q)), q\in \mA(c)\}= \{(q,c+du(q)), q\in \mA(c))\}
$$
for each $c$-critical subsolution $u$. This set is called the Aubry set, it is invariant under the flow of $H$, compact, and not empty.
It is moreover contained in the graph of the Lipschitz closed 
form $c+du$ for each $c$-critical subsolution $u$.
A consequence of the invariance of $\mA^*(c)$
is that the projected Aubry set is invariant 
under the vectorfield
$$
F(q):= \partial_pH(q,c+du(q))
$$
on $\T$
for each $c$-critical subsolution $u$.
The special $c$-critical subsolution $w$ introduced above   has the property that the strict
inequality
$$
H(q,c+dw(q))<\alpha(c) 
$$
holds on the complement of $\mA(c)$. A $c$-critical subsolution having this
property is said strict outside the Aubry set.

The function $\alpha$ is convex and superlinear on $H^1(\Tg,\R)$.
The initial definition of this function was given by John Mather in terms
of minimizing measures.
$$
\alpha(c)=\max _{\mu^*}\int_{T^*\T} (c-p)\cdot \partial_pH(q,p)+H(q,p) d\mu^*(q,p)
$$
where the maximum is taken on the set of compactly supported invariant  probability measures
$\mu^*$ on $T^*\T$.
The invariant measures minimizing this expression   will be called $c$-minimizing.
Defining the rotation number of such an invariant measure 
$$
h(\mu):= \int _{T^*\T} \partial_p H(q,p) d\mu(q,p) \subset \Rm^d=H_1(\T, \Rm),
$$
we observe that $\alpha$ is the Legendre transform 
of the function 
$$
\beta(h):= \min _{h(\mu^*)=h} \int p\cdot\partial_pH(q,p)-H(q,p) d\mu^*(q,p).
$$
In the geodesic case, where $H$ is quadratic in the fibers,
the functions $\alpha$ and $\beta$ are homogeneous of degree $2$.
The function $\sqrt{\beta}$, which is homogeneous of degree one,
is called the \textit{stable norm}.

The $c$-minimizing measures are precisely those invariant measures
which are supported on the Aubry set $\mA^*(c)$, see \cite{Man:92}.
In particular, they are supported on a Lipschitz graph, which
was a major discovery of John Mather. We denote by $\mM^*(c)$ the union 
of supports of $c$-minimizing measures. 
It is a compact invariant set contained in $\mA^*(c)$.

The subdifferential $\partial \alpha(c)$ in the sense of convex
analysis is the set of rotation numbers of $c$-minimizing measures.
When $\mu^*$ is a $c$-minimizing measure, its projection $\mu$ on $\T$
is an invariant measure of the vectorfield
$F(q)=\partial_pH(q,c+du(q))$ 
(for each $c$-critical subsolution $u$) supported on $\mA(c)$.
Its rotation number is nothing than the rotation number 
of $\mu$ as an invariant measure of $F$,
$$
h(\mu^*)=\int_{\T} F(q) d\mu(q).
$$
When $\mu^*$ is ergodic, or equivalently when $\mu$ is ergodic,
this is the asymptotic winding number of $\mu$-almost each orbit on the torus.

\subsection{The special case of dimension two, the main statement in the Tonelli case}

In this section, we work on the two-dimensional torus $\T=\Rm^2/\Zm^2$.
We recall, see \cite{FM90}, that the rotation set of a flow on the two  torus $\T$ is a compact interval contained in a straight line through the origin of $\Rm^2=H_1(\T,\Rm)$.
Moreover :
\begin{itemize}
\item If the straight line has rational direction (which means that it contains
an element of $H_1(\T,\Zm)$), then the ergodic  invariant measures 
of non-zero rotation number are supported on periodic orbits.
Moreover, the $\alpha$ and $\omega$ limit sets of the flow are made of periodic orbits.
\item If the straight line has  irrational direction, then there is at most one ergodic
invariant measure of non-zero rotation number.
\end{itemize}

Let us apply these results to the Aubry set $\mA^*(c)$
at a point $c$ which is not a minimum of the function  $\alpha$.
Then, the rotation set $\partial \alpha(c)$ does not contain 
zero, and it is contained in the rotation set of the vectorfield 
$F(q)=\partial_pH(q,c+dw(q))$ on $\T$.
We conclude that $\partial \alpha(c)$ is a compact interval of a ray
$\rho(c)\in SH_1(\T,\Rm)$, where
$$
SH_1(\T,\Rm):=(H_1(\T,\Rm)-\{0\})/]0,\infty)\approx S^1
$$ 
is the set of open  half lines of $H^1(\T,\Rm)$ starting at the origin.
We say that a ray has rational direction if it contains 
a point of $H_1(\T,\Zm)$, and
 that it has irrational direction otherwise.
\begin{itemize}
\item If $\rho(c)$ has rational direction, then the ergodic $c$-minimizing measures
are supported on periodic orbits. Moreover, the $\alpha$ and $\omega$ limits of
each orbit of $\mA^*(c)$ are periodic orbits supporting $c$-minimizing measures.
\item If $\rho(c)$ has irrational direction, then there exists a unique $c$-minimizing measure,
and the rotation set $\partial\alpha(c)$ is a point.
\end{itemize}
Observe that, in all cases, each half orbit of $\mA(c)$ has a single rotation number 
which is contained in $\rho(c)$.
However, in the case of a rational direction,
it is possible that the positive half-orbit and the negative half-orbit of a given point have different 
rotation numbers both contained in $\rho(c)$.

Let us explain a bit more how different
rotation numbers can appear in the rational case.
In this case, the periodic orbits of $\mA(c)$ are oriented  embedded 
closed curve, and they all represent the same homology class 
$[\rho(c)]\in H_1(\T,\Zm)$ which is the only indivisible integer class
in the half line $\rho(c)$.
The rotation number of the invariant measure supported on such an
 orbit is then 
$[\rho(c)]/T$, where $T$ is the minimal  period of the orbit.
The periodic orbits of $\mA(c)$ do not necessarily all have 
the same period, hence the associated measures do not 
necessarily have the same
rotation number.

For each $e>\min \alpha$, the set 
$$A(e):=\{c\in H^1(\T,\Rm): \alpha(c)\leq e\}\subset H^1(\T,\Rm)
$$
 is a compact and convex set,
whose interior is $\{\alpha < e\}$ and whose boundary is $\alpha^{-1}(e)$.
At each boundary point $c\in \alpha^{-1}(e)$,
the set $A(e)$ has a single outer normal $\rho(c)\in SH_1(\T,\Rm)$. 
The map $c\lmto \rho(c)$ is thus continuous,
and 
the  set $\alpha^{-1}(e)$ is a $C^1$ curve.
Note that the map $\rho: \alpha^{-1}(e)\lto SH_1(\T,\Rm)$
is continuous and onto, and that it preserves the order.
It is however not necessarily one to one.
For each $c\in \alpha^{-1}(e)$, we consider the face $F(c)\subset \alpha^{-1}(e)$ 
defined as the set of cohomologies $c'$ such that $\rho(c')=\rho(c)$ and $\alpha(c')=e$.
The face $F(c)$ is a compact segment containing $c$.  
It is also the set of points $c'\in A(e)$ such that 
$ (c'-c)\cdot \rho(c)=0$.
The following is well known, see \cite{Ban-90,Mas03,Mat:10}, but since we give the statement
in a way which is not obviously equivalent to those of these papers, we 
will provide a proof in section \ref{secfaces}.

\begin{prop}\label{propcases}
Each $c\in \alpha^{-1}(e)$ is in one (and only one) of the  following three cases:
\begin{enumerate}
\item \label{c1}
$\rho(c)$ has irrational direction  and $F(c)=\{c\}$.
\item \label{c2}
$\rho(c)$ has rational direction, $\mM(c)=\T$, and $F(c)=\{c\}$
\item \label{c3}
$\rho(c)$ has rational direction, $\mM(c)\neq \T$ and $F(c)$ is a non-trivial segment
$[c^-,c^+]$.
The sets $\mA^*(c^-)$ and $\mA^*(c^+)$  contain non-periodic orbits
(which are heteroclinics).
\end{enumerate}
\end{prop}

If, for a given value $e>\min \alpha$ of the energy,
case \ref{c3} does not occur for any $c\in \alpha^{-1}(e)$, then 
the map $c\lmto \rho(c)$ is a homeomorphism from 
$\alpha^{-1}(e)$ to $SH_1(\T,\Rm)$. The energy level $\{H=e\}$ is then $C^0$-integrable,
as is proved in (\cite{MS11}, Theorem 3),  see also Section  \ref{secfaces}:

\begin{prop}\label{propint}
If, for a given value $e>\min \alpha$ of the energy, 
case \ref{c3} does not occur for any $c\in \alpha^{-1}(e)$, then 
the Aubry sets $\mA^*(c), c\in \alpha^{-1}(e)$ are Lipschitz invariant graphs which 
partition the energy level $\{H=e\}$.
\end{prop}

If $H$ is the Hamiltonian associated to a Riemaniann metric, then this implies 
that the metric is flat, in view of the Theorem of Hopf, see also \cite{Innami}.
As a consequence, Theorem \ref{maingeo} follows from:

\begin{thm}\label{maintonelli}
Let $e>\min \alpha$ be a given energy level.
If there exists a cohomology $c\in \alpha^{-1}(e)$ in case \ref{c3}, then 
the polynomial entropy of the Hamiltonian flow restricted to the energy level
$\{H=e\}$ is not less than $2$.
In other words, 
if the polynomial entropy of the flow restricted to the energy level $\{H=e\}$ is 
less than two, the Aubry sets $\mA^*(c), c\in \alpha^{-1}(e)$ are Lipschitz invariant
graphs which partition the energy level.
\end{thm}

We will make use in the proof  of two important properties of the Aubry sets :

\begin{pte}\label{psemi}
The set-valed map $c\lmto \mA^*(c)$ is outer semi-continuous.
It means that  each open set $U\subset T^*\T$ containing $\mA^*(c)$,
also contains $\mA^*(c')$ for $c'$ close to $c$.
\end{pte}

\begin{pte}\label{psection}
For each $c\in \alpha^{-1}(e)$, there exists a global curve of section of 
$\mA(c)$. More precisely, there exists a cooriented $C^1$ embedded circle 
$\Sigma\subset \T$ such that $F(q)$ is transverse to $\Sigma$ on $\mA(c)$
and respects the coorientation. Moreover, each half orbit of $\mA(c)$ 
intersects $\Sigma$. The flow of $\mA(c)$ thus induces a
homeomorphism $\psi$ of $\Sigma \cap \mA(c)$ which preserves the cyclic order
of $\Sigma$. 
 \end{pte}

\textsc{Proof of Property \ref{psemi}.}
It is proved in \cite{Conley} using the content of \cite{FFR}.
\qed

\textsc{Proof of Property \ref{psection}.}
Let us consider a cohomology $c_0$ such that $\alpha(c_0)<e$ and such that
$c-c_0\in H^1(\T,\Qm)$. Note that 
$$
\rho(c)\cdot (c-c_0)>0.
$$
We consider a $c$-critical subsolution $w$ of class $C^{1,1}$ and 
strict outside the Aubry set.
We also consider a $c_0$-critical subsolution $u_0$.
Let $l\in \Nm$ be such that $l(c-c_0)\in H^1(\T,\Zm)$.
Let us consider the $C^{1,1}$ function $\hat \Theta$ on $\Rm^2$ defined by
$$
q\lmto \hat \Theta(q)=l(c-c_0)q+lw(q)-lu_0(q). 
$$
The function $\tilde \Theta =\hat \Theta \mod 1 : \Rm^2\lto \Tm=\Rm/\Zm$
is $\Zm^2$-periodic, hence it gives rise to a function 
$\Theta: \T\lto \Tm$ such that 
$d\Theta=l(c-c_0)+l(dw-du_0)$.
Let us use as above the notation $F(q)=\partial_pH(q,c+dw(q))$.
For each point $q$ such that $H(q,c+dw(q))=e$, we have
\begin{align*}
e-l^{-1}d\Theta(q)\cdot F(q) &=H(q,c+dw(q))- \partial_pH(q,c+dw(q))\cdot (c_0+du_0(q)-c-dw(q))\\
& \leq H(q,c_0+du_0(q))\leq \alpha(c_0)<e
\end{align*}
hence 
$$
d\Theta (q) \cdot F(q) >0.
$$
Let us consider a regular value $\theta$ of $\Theta$.
Such a value exists by fine versions of  Sard's Theorem (see \cite{Ba:93}) since $\Theta$
is $C^{1,1}$.
The preimage $\Theta^{-1}(\theta)$ is a $1$-dimensional cooriented
submanifold of $\T$. It can be seen as an intersection cocycle of cohomology $l(c-c_0)$.
It is a finite union of embedded cooriented circles 
$\Sigma_i$ each of which is a cocycle of cohomology $\sigma_i$, with $\sum \sigma_i=l(c-c_0)$.
Since $\rho(c) \cdot l(c-c_0)>0$, there exists $j$ such that 
$\rho(c)\cdot \sigma_j>0$.
We denote by $\Sigma$ the cooriented circle $\Sigma_j$.
We have $d\Theta (q) \cdot F(q) >0$ hence the orbits of $\mA(c)$ 
are transverse to $\Sigma$, and intersect it according to the coorientation.

Finally, each half orbit of $\mA(c)$ has a rotation number 
contained in $\rho(c)$.
 We have seen that $\sigma\cdot \rho(c)>0$, where $\sigma$ is the cohomology 
of the intersection cocycle associated to $\Sigma$. 
Each half orbit of $\mA(c)$
thus intersects $\Sigma$.
As a consequence, the flow of $\mA(c)$  generates a Poincar\'e map 
$$
\psi: \mA(c)\cap \Sigma\lto \mA(c)\cap \Sigma
$$
which is a bi-Lipschitz homeomorphism preserving the cyclic order on the circle $\Sigma$.
This implies that $\psi $ can be extended to a homeomorphism of $\Sigma$
preserving the cyclic order.
\qed

\subsection{Faces of the balls of $\alpha$ on the two torus.}\label{secfaces}

We take $d=2$ and fix an energy level $e>\min \alpha$.
We study the affine parts of the ball $\alpha^{-1}(e)$ and 
prove Propositions \ref{propcases} and \ref{propint}.
The following is a variant of  a Lemma of Daniel Massart \cite{Mas03}:

\begin{lem}\label{lemmatherface}
Assume that $d=2$, that $e>\min \alpha$, and that  
the  segment $[c_0,c_1]$ is contained in $\alpha^{-1}(e)$.
Then the Mather set $\mM^*(c)$ does not change when $c$ varies in 
$[c_0,c_1]$. 
\end{lem}

\proof
We have $(c_1-c_0)\cdot \rho(c)=0$ for each $c\in [c_0,c_1]$.
Let $c$   be a point in $[c_0,c_1]$ and let
$\mu^*$ be an ergodic $c$-minimizing measure.
Such a measure has a rotation number
$h(\mu^*)=s\rho(c), s>0$.
For each $c'$, we have
$$
\alpha(c')\geq 
\int_{T^*\T} (c'-p)\cdot \partial_pH(q,p)+H(q,p) d\mu^*(q,p)
=\alpha(c)+(c'-c)\cdot h(\mu^*).
$$
If $c'\in [c_0,c_1]$, then the inequality 
$\alpha(c')\geq \alpha(c)+(c'-c)\cdot h(\mu^*)$
is an equality, 
hence $\mu^*$ is a $c'$-minimizing measure.
This implies that $\mM^*(c')\subset \mM^*(c)$ for each $c'\in [c_0,c_1]$.
By symetry, we conclude that 
$\mM^*(c')=\mM^*(c)$ for $c'\in [c_0,c_1]$.
\qed

The following Lemma also comes from Daniel Massart \cite{Mas03}:

\begin{lem}\label{lemaubryface}
Assume  that $e>\min \alpha$, and that the  segment $[c_0,c_1]$ is contained in $\alpha^{-1}(e)$.
Then the Aubry set $\mA^*(c)$ does not change when $c$ varies in $]c_0,c_1[$. Moreover, we have the inclusion  
$\mA^*(c)\subset \mA^*(c_0) \cap \mA^*(c_1)$ for each $c\in\, ]c_0,c_1[$.
\end{lem}

We will see that, unlike the Mather set, the Aubry set can be bigger at the boundary points.

\proof
Let us consider a point $c=ac_0+(1-a)c_1, a\in ]0,1[$.
Let $w_i$,${i\in \{0,1\}}$ be a $c_i$-critical subsolution strict outside the Aubry set.
Then, $w_c:=aw_0+(1-a)w_1$ is a $c$-critical subsolution.
Using the strict convexity of $H$ in $p$, we observe that the strict inequality $H(q,c+sw_c(q))<e$ holds outside of the set where $H(q,c_0+dw_0(q))=e$ and  $H(q,c_1+dw_1(q))=e$ and $c_0+dw_0= c_1+dw_1$.
We conclude that the Aubry set $\mA(c)$ is contained in $\mA(c_0)\cap \mA(c_1)$, and that $c_0+dw_0=c_1+dw_1=c+dw_c$ on $\mA(c)$.
As a consequence, $\mA^*(c)\subset \mA^*(c_0)\cap\mA^*(c_1)$.
If $c$ and $c'$ are two points in $]c_0,c_1[$, assuming for example that $c\in \,]c_0,c'[$, we conclude that $\mA^*(c)\subset \mA(c')$. 
Similarly, we have $c'\in\, ]c,c_1[$, and we obtain the converse inclusion,
hence $\mA^*(c)=\mA^*(c')$.
\qed

We recall that $F(c)$ is defined as the largest segment of $\alpha^{-1}(e)$ containing $c$.

\begin{cor}\label{corcase2}
If $\mM(c)=\T$, then $F(c)=c$.
\end{cor}

\proof
If $\mM(c)=\T$, then there exist one and only one $c$-critical subolution $w$,
and $\mM^*(c)$ is the graph of $c+dw$.
Assume now that there exists $c'$ such that $\mM^*(c')=\mM^*(c)$.
Then $\mM(c')=\T$, hence there exists a unique $c'$-critical subsolution $w'$, 
and $\mM^*(c')$ is the graph of $c'+dw'$. We thus have $c+dw=c'+dw'$ at each point,
which implies that $c=c'$.
In view of  Lemma \ref{lemmatherface}, this implies that $F(c)=c$. 
\qed

The following was first proved by Bangert, see \cite{Ban-90}:

\begin{cor}\label{corcase1}
If $\rho(c)$ has an irrational direction, then $F(c)=c$.
\end{cor}

\proof
As above, let us consider a cohomology $c$ satisfying the hypothesis
of the Corollary, and a cohomology $c'$ such that
$[c,c']\in \alpha^{-1}(e)$, hence $\mM^*(c')=\mM^*(c)$, by Lemma \ref{lemmatherface}.
We will prove that $c'=c$, which implies the Corollary.
Note that $(c'-c)\cdot \rho(c)=0$, so that
it is enough to prove that $(c'-c)\cdot[\Sigma]=0$,
where $[\Sigma]\in H_1(\T,\Zm)$ is the homology of the section $\Sigma$
given by \ref{psection} (equipped with an orientation).
Let $w$ and $w'$ be $c$ and $c'$-critical subsolutions.
We consider the closed Lipschitz form $\omega=c'-c+dw'-dw$, 
whose cohomology is $c'-c$ and prove that $\int_{\Sigma}\omega =0$.

Since $\mM^*(c')=\mM^*(c)$, we have $\omega=0$ on $\mM(c)$.
It is thus enough to prove that 
 $\int_{I}\omega =0$ 
for each connected component  $I$
 of the complement of $\mM(c)\cap \Sigma$ in 
$\Sigma$.
We first observe that 
$$
\int_{I}\omega =\int _{\psi(I)} \omega,
$$
where $\psi$ is a homeomorphism of $\Sigma$ extending the return map 
of $\Sigma \cap \mA(c)$.
To prove this equality, we integrate $\omega$ on  the contractible closed curve made of the interval $I=\,]q^-,q^+[$, followed by the orbit of $q^+$ until its return $\psi(q^+)$, followed by the interval $-\psi(I)=\,]\psi(q^+), \psi(q^-)[$
followed by the piece of orbit of $q^-$ in negative time direction from $\psi(q^-)$ to $q^+$.

Since the intervals $\psi^k(I)$ are two by two disjoint in $\Sigma$, 
their lengh is converging to zero.
Since the form $\omega$ is bounded, this implies that 
$\int_{\psi^k(I)}\omega \lto 0$, hence that $\int_I \omega=0$. 
\qed

In view of these corollaries there are three cases:
\begin{itemize}
\item
$\rho(c)$ has irrational direction and $F(c)=\{c\}$ 
(Corollary \ref{corcase1}).
\item
$\rho(c)$ has rational direction, $\mM(c)=\T$, and 
$F(c)=\{c\}$ 
(Corollary \ref{corcase2}).
\item
$\rho(c)$ has rational direction and $\mM(c)\neq \T$.
\end{itemize}

Let us study more precisely the last case.
We denote by $[c^-,c^+]$ the face $F(c)$.
We assume for definiteness that $c$ is an interior point of this face, which means that either $c\in\, ]c^-,c^+[$ or $c^-=c=c^+$.

We consider the cooriented section $\Sigma$
given by Property \ref{psection}.
We  orient $\Sigma$ in such a way that $(c^+-c^-)\cdot [\Sigma]\geq 0$,
where $[\Sigma]$ is the homology of $\Sigma$
(hence  $(c^+-c^-)\cdot [\Sigma]> 0$ if $c^+\neq c ^-$, since
$[\Sigma]$ is not proportional to $\rho(c)$).
The return map $\psi$ from $\mM(c)\cap \Sigma$ to istelf
is periodic, its minimal period is $\tau:=\sigma \cdot [\rho(c)]$.
The complement of $\mM(c)$ in $\T$ is a disjoint union of topological open annuli.
Each of these annuli $U$ intersects $\Sigma$ in $\tau$ disjoint open intervals, that we orient
according to the orientation of $\Sigma$.
Each orbit of $\mA(c)-\mM(c)$ is contained in such an annulus $U$,
is $\alpha$-asymptotic to one of its boundaries, and is $\omega$-asymptotic to
its other boundary. 
We say that such an orbit is positive if it crosses the annulus $U$ according to the orientation of 
$\Sigma$, and that it is negative if it crosses in the other direction.
In other words, the heteroclinic orbit is positive if the sequence of its successive intersections
with the interval $I$ is increasing.
The following implies Proposition \ref{propcases}:

\begin{prop}\label{propheteroclines}
If $c\in \alpha^{-1}(e)$ is such that $\rho(c)$ is rational and $\mM(c)\neq \T$, then 
\begin{itemize}
\item
$c^-\neq c^+$ 
\item
The Aubry set $\mA(c^+)$ contains positive heteroclinics in each 
connected component of $\T-\mM(c)$ and no other orbit except those of $\mM(c)$.
\item
The Aubry set $\mA(c^-)$ contains negative heteroclinics in each 
connected component of $\T-\mM(c)$ and no other orbit except those of $\mM(c)$.
\item
Finally, $\mA(c)=\mM(c)$ for each $c$ in $]c^-,c^+[$.
\end{itemize}
\end{prop}

\proof
The statements concerning $\mA(c^+)$ and $\mA(c^-)$ imply that 
$c^+\neq c^-$.
Moreover, they imply that
$\mA(c^+)\cap \mA(c^-)=\mM(c)$, hence that $ \mA(c)=\mM(c)$ for each $c$ in $]c^-,c^+[$,
by Lemma \ref{lemaubryface}.

We will now prove the statement concerning $c^+$, the other one being similar.
We fix a connected component $U$ of $\T-\mM(c)$, and a connected component
$I$ of $U\cap \Sigma$.
Let $\rho_m$ be the direction of $m[\rho]+[\Sigma]$, and let $c_m\in \alpha^{-1}(e)$ be such that 
$\rho(c_m)=\rho_m$. Note that $c_m\lto c^+$.

The annulus $U$ contains an oriented closed curve $K$ of homology $[\rho]$.
The Aubry set $\mA(c_m)$ contains a periodic orbit of homology positively proportional to $m[\rho]+[\Sigma]$,
hence it intersects $K$.
In view of the semi-continuity of the Aubry set, see Property \ref{psemi}, we deduce that $\mA(c^+)$
intersects $K$. 
As a consequence, the set $\mA(c^+)$ does contain heteroclinic orbits 
in $U$.
Moreover, there exists $m_0\in \Nm$ and  a compact subinterval $J\subset I$
which contains a point in each orbit  of $\mM(c_m)$ for $m\geq m_0$.

Let us now consider the return map $\psi^{\tau}$ of  $I\cap \mA(c)$.
Either we have $\psi^{\tau}(x)> x$ for each $x\in I\cap \mA(c)$, and the orbits of $\mA(c)\cap U$
are positive heteroclinics,
or we have $\psi^{\tau}(x)< x$ for each $x\in I\cap \mA(c)$, and the orbits of $\mA(c)\cap U$
are negative heteroclinics.

In view of the semi-continuity of the Aubry set, $\Sigma$ is also a cooriented transverse section for 
$\mA(c_m)$ for $m\geq m_0$, provided $m_0$ is large enough.
 Denoting by $\psi_{c_m}$ the corresponding section map,
 we have $\psi^{\tau}_{c_m}(J)\subset I$ for $m\geq m_0$ provided $m_0$ is large enough.
Then, there  exists a point $x_m\in J\cap \mM(c_m)\subset J\cap \mA(c_m)$, 
and $x_m<\psi^{\tau}_{c_m}(x_m)$.

At the limit, using the semi-continuity of the Aubry set, we find a point $x\in J\cap \mA(c^+)$ such that
$\psi^{\tau}(x)\geq x$, hence $\psi^\tau(x)>x$. We conclude that the heteroclinics are positive.
\qed

For the convenience of the reader,  and because our statement is not exactly the one of \cite{MS11}
Theorem 3, we now prove Proposition \ref{propint}, following  \cite{MS11}:

We consider an energy level $e> \min \alpha$ such that the curve $\alpha^{-1}(e)$ does not contain any non-trivial segment,
which is equivalent to saying that $\mM(c)=\T$ for each $c$ such that $\rho(c)$ is rational.
Note then that the map $\rho: \alpha^{-1}(e)\lto SH_1(\T,\Rm)$ is continuous and bijective, hence it is  a homeomorphism.
Since the set $SH_1(\T,\Zm)$ of rational directions is dense in $SH_1(\T,\Rm)$, its preimage $\rho^{-1}(SH_1(\T,\Zm))$ is dense in $\alpha^{-1}(e)$.
For each point $c$ in this set, we have $\mA(c)=\T$.
In view of the semi-continuity of the Aubry set, we deduce that $\mA(c)=\T$ for each $c\in \alpha^{-1}(e)$.
As a consequence, there exists a unique (up to the addition of a constant) $c$-critical subsolution $w_c$, which is actually a solution, and  the Aubry set $\mA^*(c)$ is the graph of $c+dw_c$. 
Moreover, the functions $dw_c, c\in \alpha^{-1}(e)$ are equi-Lipschitz.
The semicontinuity of the Aubry set $\mA^*$ implies that the map $c\lmto c+dw_c(q)$ is continuous for each $q\in \T$.

The orbits of $\mA^*(c)$ all have a forward rotation number in $\rho(c)$.
For $c'\neq c$, the orbits of $\mA^*(c')$ all have a forward rotation number in $\rho(c')$,
and, since $\rho(c')\neq \rho(c)$, the sets $\mA^*(c)$ and $\mA^*(c')$ are disjoint.
As a consequence, for each $q\in \T$, the map $c\lmto c+dw_c(q)$ is one to one on  $\alpha^{-1}(e)$, hence it has degree $\pm 1$ as a circle map into $\{p\in T_q\T : H(q,p)=e\}$.
It is thus onto, which implies that the Aubry sets fill the energy level.
\qed

\section{Lower bound for the polynomial entropy}\label{secproof}

We prove Theorem \ref{maintonelli}.
We consider an energy level $e>\min \alpha$,  assume that the ball $\alpha^{-1}(e)$ contains
a non-trivial face $[c^-,c^+]$, and prove that the entropy of the Hamiltonian flow on the energy level 
$H^{-1}(e)$ is at least two. The proof have similarities with the ones of \cite{Mar13}, \cite{LM}, and \cite{Mar14}.
We work with the section $\Sigma$ of $\mA(c^+)$ given in Property \ref{psection}.
We fix a parameterisation $\Rm/\Zm\lto \Sigma$, and put on $\Sigma$
the distance such that this parameterisation is  isometric.
This distance is Lipschitz equivalent to the restriction
of the distance on $\T$.

The direction $\rho(c)\in SH_1(\T,\Zm)$ is independant of $c\in [c^-,c^+]$,
and it is rational, we denote it by $\rho$ in the sequel, and by 
$[\rho]\in H_1(\T,\Zm)$ the associated indivisible integer point.
As above, we consider, for $m\in \Nm$, the direction $\rho_m$ of $m[\rho]+[\Sigma]$ and a cohomology $c_m\in \alpha^{-1}(e)$
such that $\rho(c_m)=\rho_m$. Let us decide for definiteness that $c_m\lto c^+$ (otherwise we exchange the names of $c^+$
and $c^-$).
We fix $m_0$ large enough so that $\Sigma$ is a cooriented transverse section of $\mA(c_m)$
for each $m\geq m_0$ (such a value of $m_0$ exists in view of the semi-continuity of the Auby set)
and denote by $\psi_{c_m}$ the corresponding return map of $\mA(c_m)\cap \Sigma$.
The orbits of $\mM(c)$ give rise to periodic orbits of $\psi$, and
   the minimal period of these orbits is
 $\tau:= \sigma \cdot [\rho]$, where $\sigma$ is the cohomology of
 the intersection cocycle associated to $\Sigma$.
 The integral class $m[\rho]+[\Sigma]$ is not necessarily indivisible in $H_1(\T,\Zm)$,
 hence the minimal period for the return map $\psi_{c_m}$ of the points of $\Sigma \cap \mM(c_m)$ 
 may be smaller than $m\tau$.
 However,  because $m[\rho]+[\Sigma]$ is indivisible in the group generated by $[\Sigma]$ and $[\rho]$,
 we have:
 
 \begin{lem}
 The orbits of $\mM(c_m)$ give rise to periodic orbits of $\psi^{\tau}_{c_m}$ with minimal period $m$.
 \end{lem}
 
 The following proof might appear unnecessarily complicated. 
 Things can also be understood  as follows:
 The statement of the Lemma is obvious if ($[\Sigma]$,$[\rho])$ is a base of $H^1(\T,\Zm)$
 (since $m[\rho]+[\Sigma]$ is then indivisible in $H^1(\T,\Zm)$)
 and we can reduce the situation to this simple case by taking a finite
 covering, which does not change the value of the polynomial entropy.

 \proof
Let us denote by $G$ the subgroup of $H_1(\T,\Zm)$ generated by $[\Sigma]$ and $[\rho]$.
Let $\chi: \T\lto \T$ be a covering such that $\chi_*(H_1(\T,\Zm))=G$. This coverging has $\tau$ sheets.
The preimage of $\Sigma$ by this coverging has $\tau$ connected components, and we denote by $\tilde \Sigma$ one of them. 
We have 
 $$
 [\tilde \Sigma]=\chi_*^{-1}([\Sigma]) \subset H_1(\T,\Zm),
 $$
and the cohomology of the intersection cocycle associated to $\tilde \Sigma$ is $\tilde \sigma := \chi^*(\sigma)/\tau$.
The preimage of each orbit of $\mM(\rho)$ has $\tau$ connected components, each of which is a closed curve of homology 
$[\tilde \rho]:= \chi_*^{-1}([\rho])$.
Since $[\rho_m]$ does not necessarily belong to $G$, the preimage of closed orbits in $\mM(c_m)$ may have less than $\tau$
connected components. 
Each of these components have a homology in  $H_1(\T,\Zm)$ which is indivisible and positively proportional to 
$
[\tilde \rho_m]:= m[\tilde \Sigma]+[\tilde \rho].
$
Since $m[\rho]+[\Sigma]$ is indivisible in $G$, $[\tilde \rho_m]$ is indivisible in $H_1(\T,\Zm)$, hence it is equal to the homology of the connected components of the preimages of orbits of $\mM(c_m)$.
Note that $\chi_*([\tilde \rho_m])=m[\rho]+[\Sigma]$ is not necessarily equal to $[\rho_m]$.

The minimal $\psi_{c_m}^{\tau}$-period of orbits of $\mM(c_m)\cap \Sigma$ is equal to the intersection number 
$\tilde \sigma \cdot [\tilde\rho_m]=(\sigma/\tau) \cdot (m[\rho]+[\Sigma])=m$.
 \qed

As in the proof of Proposition \ref{propheteroclines},
we consider a compact subinterval  $J\subset I$  such that 
each orbit of $\mA(c_m)$ contains a point of $J$ for $m\geq m_0$
(we may have to increase $m_0$).
We can chose $\epsilon_0>0$ small enough  and $m_0$
large enough to have, for all $m\geq m_0$,
\begin{equation}\label{epsilon0}
d(q,\psi^{\pm \tau}_{c_m}(q))\geq 2\epsilon_0 
\end{equation}
for all
$q\in \mA(c_m)\cap I$ such that $d(q,J)\leq \epsilon_0$.
Each point of  $J\cap \mM(c_m), m\geq m_0$  is at distance at least $2\epsilon_0$
 from all the other  points of its $\psi^{\tau}$-orbit.
We deduce: 

\begin{lem}\label{orbit}
For $q\in J\cap  \mM(c_m), m\geq m_0$, the orbit 
$$
O_{\psi^{\tau}}(q):=
 \{
q, \psi^{\tau}_{c_m}(q), \psi^{2\tau}_{c_m}(q),\ldots, \psi^{(m-1)\tau}_{c_m}(q)
\}
$$
is $(\epsilon_0,m)$-separated by $\psi^{\tau}$, hence 
 $(2\epsilon_0,\tau m)$-separated by $\psi$. 
\end{lem}

\proof
Let $\theta$ and $\theta'$ be two points of this orbit.
There exists $l\in \{0,1,\ldots,m-1\}$ such that $\psi^{l\tau}(\theta)\in J$. Then 
$\psi^{l\tau}(\theta')$ is another point of the same orbit, hence 
$d\big(\psi^{l\tau}(\theta),\psi^{l\tau}(\theta')\big)\geq 2\epsilon_0$.
\qed

We denote by $\Sigma^*$ the set of points of the energy surface $H^{-1}(e)$ which project on $\Sigma$ ,
we endow it with a distance which satisfies $d((q,p),(q',p'))\geq d(q,q')$.
We consider the compact invariant set of the Hamiltonian flow (on the energy level) defined by
$$
A:= \mA^*(c^+)\cup\bigcup_{m\geq m_0} \mA^*(c_m).
$$
The surface  $\Sigma^*$ is a transverse section of the flow on this invariant set, as required  in 
Proposition \ref{Poincare}. We denote by $\Psi$ the corresponding return map of $A\cap \Sigma^*$.
The restriction of $\Psi$ to $\mA^*(c_m)\cap \Sigma^*$ is conjugated to $\psi_{c_m}$ by the projection.
In view of Proposition \ref{Poincare}, it is enough to bound from below the polynomial entropy of $\Psi$
on $A\cap \Sigma^*$. We  exhibit a sufficiently large separated set using the orbits
$$
O_{\Psi^{\tau}}(x):=\{x, \Psi^{\tau}(x), \ldots, \Psi^{(k-1)\tau}(x)\}.
$$

\begin{lem}
Let us chose, for each $m\geq m_0$,  an element $x_m$ of $\mM^*(c_m)\cap J^*$.
The set 
$$
\bigcup_{k\in \{m,m+1,m+2,\ldots, 2m-1\}} O_{\Psi^{\tau}}(x_k)
$$
is $(\epsilon_0,4m)$-separated by $\Psi^{\tau}$, hence $(\epsilon_0,4m\tau)$-separated by $\Psi$ (for $m\geq m_0$).
\end{lem}

Since the cardinal of this union is more than $m^2$, we conclude that the polynomial entropy of $\Psi$ is at least two.
By Proposition \ref{Poincare} the polynomial entropy of the Hamiltoinan flow on the energy surface is at least two.

\proof
Let $x=(q,p)\in O_{\Psi^{\tau}}(x_k)$ and $y=(\theta, \eta)\in O_{\Psi^{\tau}}(x_l)$ be two different points in this union.

If $k=l$, the points  $x$ and $y$ belong to the same orbit $O_{\Psi^{\tau}}(x_k)$.  They are 
$(\epsilon_0,m)$-separated by $\Psi^{\tau}$ in view of Lemma \ref{orbit}.

Otherwise, we assume for definiteness that $m\leq k<l< 2m$.
There exists an integer $s\in \{0,1,\ldots, k-1\}$ 
such that
$\psi^{s\tau}_{c_k}(q)\in J$.

If $d\big(\psi^{s\tau}_{c_k}(q),\psi^{s\tau}_{c_l}(q)\big)\geq \epsilon_0$, then 
$d\big(\Psi^{s\tau}(x), \Psi^{s\tau}(y)\big)\geq \epsilon_0$ hence 
$x$ and $y$ are $(\epsilon_0, k)$-separated by $\Psi^{\tau}$.

If  $d\big(\psi^{s\tau}_{c_k}(q),\psi^{s\tau}_{c_l}(\theta)\big)\leq \epsilon_0$, then 
$d\big(\psi^{s\tau}_{c_l}(\theta),J\big)\leq \epsilon_0$  hence
 (\ref{epsilon0}) implies that
$$
d\big( \psi^{(k+s)\tau}_{c_l}(\theta) ,\psi^{s\tau}_{c_l}(\theta)
\big)\geq 2\epsilon_0
$$
so that
$$
d\big( \psi^{(k+s)\tau}_{c_l}(\theta) ,\psi^{(k+s)\tau}_{c_k}(q)
\big)=
d\big( \psi^{(k+s)\tau}_{c_l}(\theta) ,\psi^{s\tau}_{c_k}(q)
\big)
\geq \epsilon_0
$$
hence
$$
d\big( \Psi^{(k+s)\tau}(y) ,\Psi^{(k+s)\tau}(x)
\big)\geq \epsilon_0.
$$
As a consequence, the points $x$ and $y$ are $(\epsilon_0,4m)$-separated by $\Psi^{\tau}$.
\qed

\end{document}